\magnification=\magstep1
\input amstex
\documentstyle{amsppt}
\pagewidth{6.5truein}
\pageheight{9truein}
\loadbold

\def\R{{\bold R}}
\def\C{{\bold C}}
\def\N{{\bold N}}

\def\id{\mathop{\roman{id}}\nolimits}

\def\Im{\mathop{\roman{Im}}\nolimits}

\def\graph{\mathop{\roman{graph}}\nolimits}

\topmatter
\title Analytic and Nash equivalence relations of Nash maps\endtitle 
\author Masahiro SHIOTA \endauthor 
\address 
Graduate School of Mathematics, Nagoya University, Chikusa, Nagoya, 
464-8602, Japan\endaddress
\abstract
Let $M$ and $N$ be Nash manifolds, and $f$ and $g$ Nash maps from $M$ to $N$. 
If $M$ and $N$ are compact and if $f$ and $g$ are analytically R-L equivalent, then they are Nash R-L equivalent. 
In the local case, $C^\infty$ R-L equivalence of two Nash map germs implies Nash R-L equivalence. 
This shows a difference of Nash map germs and analytic map germs. 
Indeed, there are two analytic map germs from $(\R^2,0)$ to $(\R^4,0)$ which are $C^\infty$ R-L equivalent but not 
analytically R-L equivalent. 
\endabstract 
\email
shiota\@math.nagoya-u.ac.jp
\endemail
\endtopmatter

\document
\head \S 1. Introduction
\endhead
A {\it Nash manifold} is a semialgebraic and analytic submanifold of a Euclidean space. 
A {\it Nash map} between Nash manifolds is an analytic map with semialgebraic graph. 
A {\it Nash set} in or a {\it Nash subset} of a Nash manifold is the zero point set of a Nash function on the 
manifold. 
A {\it Nash closure} of a subset of a Nash manifold $M$ ({\it in} $M$) is the smallest Nash set in $M$ 
containing the subset. 
Let $f$ and $g$ be Nash maps from a Nash manifold $M_1$ to another $M_2$. 
We say that $f$ and $g$ are {\it Nash R-L equivalent} if there exist Nash diffeomorphisms $\tau_1$ of $M_1$ and 
$\tau_2$ of $M_2$ such that $f\circ\tau_1=\tau_2\circ g$. 
If $\tau_2=\id$, we say that $f$ and $g$ are {\it Nash R equivalent}. 
In the same way we define {\it analytic $(C^\infty)$ R-L} and {\it R equivalence} of two analytic ($C^\infty$) 
maps between analytic ($C^\infty$, resp.) manifolds. \par
  Classification of maps (map germs) by R-L equivalence relation seems to be natural and is more difficult to solve than one by R equivalence relation. 
(See [T$_1$]). 
Specialists of real singularity theory state theorems about $C^\infty$ or analytic maps (map germs), but consider in mind Nash or polynomial maps (map germs) except when they show pathological phenomena. 
Moreover, we do not expect good theory on classification of polynomial maps (map germs) by polynomial R-L equivalence relation. 
Hence it is worth constructing theory of classification of Nash maps (map germs) by Nash R-L equivalence relation. 
Then we need to avoid integration of vector fields. 
Historically, this has been one of the most useful methods of classification by $C^\infty$, analytic or topological R-L equivalence relation, e.g., Mather's work on stability of $C^\infty$ maps and [T$_{1,2}$]. 
However, the integration of a Nash vector field is very seldom of class Nash. 
Hence we want to know whether $C^\infty$ or analytic R-L equivalence of two Nash maps (map germs) implies Nash R-L equivalence. 
The main theorem is the following. 
\proclaim{Theorem} 
Let $f$ and $g$ be analytically R-L equivalent Nash maps between Nash manifolds. 
If the manifolds are compact, $f$ and $g$ are Nash R-L equivalent. 
\endproclaim
  The local case for R equivalence immediately follows from Artin Approximation Theorem [A$_2$], which is called {\it AA Theorem}. 
To be precise, for Nash map germs $\phi,\psi:(\R^n,0)\to(\R^m,0)$, if there exists an analytic diffeomorphism 
germ $\tau$ at 0 in $\R^n$ such that $\phi\circ\tau=\psi$ then $\tau$ is approximated by a Nash diffeomorphism 
germ $\tilde\tau$ in the $\frak p$-adic topology so that $\phi\circ\tilde\tau=\psi$, where $:(\R^n,0)\to(\R^m,0)$ 
denotes the germ at 0 of a map from a neighborhood of 0 in $\R^n$ to one in $\R^m$ carrying 0 to 0 and $\frak p$ denotes the maximal ideal of the ring of convergent power series in $n$-variables. 
The global case for R equivalence was shown in [C-R-S], where AA Theorem is globalized to G (=Global) AA Theorem by 
N\'eron Desingularization Theorem. 
For the proof of the above theorem we use a ``nested type'' of GAA Theorem in [F-S$_2$], which is called {\it NGAA Theorem}. \par
  In the theorem we cannot replace the assumption that $f$ and $g$ are analytically R-L equivalent by the one 
that they are $C^\infty$ R-L equivalent as follows. 
Let $f:S^1\to S^2$ be a Nash map such that $\Im f$ is a simple Nash curve and its Nash closure in $S^2$ is Nash 
diffeomorphic to $S^1$. 
Let $X$ be a smooth simple Nash curve in $S^2$ whose Nash closure is not Nash diffeomorphic to $S^1$ (e.g., a 
$C^\infty$ smooth simple subcurve of an algebraic curve in $S^2$ which is homeomorphic but not $C^\infty$ 
diffeomorphic to $S^1$), and $\pi:\Im f\to X$ a Nash diffeomorphism, which exists by Theorem VI.2.2 in [S$_1$]. 
Set $g=\pi\circ f$. 
Then $f$ and $g$ are $C^\infty$ R-L equivalent because $\pi$ is extended to a $C^\infty$ diffeomorphism of $S^2$, 
but they are not Nash R-L equivalent because there does not exist a Nash diffeomorphism of $S^2$ which carries 
$\Im f$ to $X$. 
This phenomenon happens because we require the diffeomorphism of $S^2$ to be globally of class Nash. 
In the local case, the assumption of $C^\infty$ R-L equivalence is sufficient (theorem 4).\par
  The compactness assumption in the theorem is also necessary. 
Indeed, there exist two polynomial functions on $\R^8$ which are analytically R equivalent but not Nash R-L equivalent 
(II.7.13 in [S$_2$]). \par
  In section 2 we prove the theorem. 
In section 3 we treat the local cases and show the two facts stated in the abstract. 
  See [S$_3$] and [F-S$_1$] for other results on equivalence of maps (map germs) by similar view points and for the real analytic and Nash sheaf theory, which we need in the proof below. \par
  Two $C^\infty$ R-L equivalent $C^\omega$ map germs are not necessarily $C^\omega$ R-L equivalent but two $C^\infty$ R-L equivalent Nash map germs are always Nash R-L equivalent. 
We can say the reason is that the image of an analytic set germ under an analytic map germ is not necessarily semianalytic but the image of a Nash set under a Nash map is semialgebraic. 

\head \S 2. Proof of the theorem
\endhead
  We prove the theorem in a more general form. 
Let $M_1\supset X_1,\,M_2\supset X_2,\,L_1\supset Y_1$ and $L_2\supset Y_2$ be Nash manifolds and closed 
semialgebraic subsets. 
Let $N(M_1)$ denote the topological space of Nash functions on $M_1$ with the compact-open $C^\infty$ topology, 
$\Cal N_{M_1}$ the sheaf of Nash function germs at points in $M_1$, $\Cal N(X_1)$ the topological space of 
germs on $X_1$ of Nash functions defined on semialgebraic neighborhoods of $X_1$ in $M_1$ with the topology of 
the inductive limit space of $N(U)$ where $U$ runs through the family of open semialgebraic neighborhoods of 
$X_1$ in $M_1$, and $\Cal N(X_1,X_2)$ the topological space of germs on $X_1$ of Nash maps from semialgebraic 
neighborhoods of $X_1$ in $M_1$ to ones of $X_2$ in $M_2$ which carry $X_1$ to $X_2$ with the topology defined 
in the same way. 
Note that $\Cal N(X_1)$ is a metrizable linear topological space, and when we regard $M_1$ as a Euclidean 
space by a local coordinate neighborhood, a sequence $f_n,\,n=1,2,...,$ in $\Cal N(X_1)$ converges to 0 if and 
only if for each derivative $D^\alpha$ and for each compact subset $K$ of $X_1$, the sequence of the 
restrictions to $K$ of $D^\alpha f_n$ converges uniformly to 0. 
A {\it Nash diffeomorphism germ} in $\Cal N(X_1,Y_1)$ is the germ on $X_1$ of a Nash diffeomorphism $\tau_1$ 
from a semialgebraic neighborhood of $X_1$ in $M_1$ to one of $Y_1$ in $L_1$ such that $\tau_1(X_1)=Y_1$. 
Let $f\in\Cal N(X_1,X_2)$ and $g\in\Cal N(Y_1,Y_2)$. 
We say that $f$ and $g$ are {\it Nash R-L equivalent} if there exist Nash diffeomorphism germs $\tau_1\in
\Cal N(X_1,Y_1)$ and $\tau_2\in\Cal N(X_2,Y_2)$ such that $\tau_2\circ f=g\circ\tau_1$. 
In the same way we define $\Cal O_{M_1}$ the sheaf of analytic function germs at points in $M_1$, $\Cal O(X_1
)$ the topological space of analytic function germs on $X_1$, $\Cal O(X_1,X_2)$ the topological space of 
analytic map germs on $X_1$ and the {\it analytic R-L equivalence}. 
Then we generalize the theorem as follows. 

\proclaim{Theorem 1} 
Let $f\in\Cal N(X_1,X_2)$ and $g\in\Cal N(Y_1,Y_2)$. 
Assume that $X_1,\,X_2,\,Y_1$ and $Y_2$ are compact. 
If $f$ and $g$ are analytically R-L equivalent, then they are Nash R-L equivalent. 
\endproclaim
The case of R equivalence is Theorem 3.2 in [F-S$_2$]. 
We proceed to prove in the same way as in [F-S$_2$]. 
The key is the following NGAA Theorem (Proposition 3.1 in [F-S$_2$]). 

\proclaim{Theorem 2} 
Let $M_1,...,M_m$ be Nash manifolds, $X_1\subset M_1,...,X_m\subset M_m$ compact semialgebraic subsets, and 
$l_1,...,l_m\in\N\,(=\{0,1,...\})$. 
Let $F_i(x_1,...,x_i,y_1,...,y_i)\allowmathbreak\in\Cal N(X_1\times\cdots\times X_i\times\R^{l_1}\times\cdots
\times\R^{l_i})$ and $f_i(x_1,...,x_i)\in\Cal O(X_1\times\cdots\times X_i)^{l_i}$ for $i=1,...,m$ such that 
$F_i(x_1,...,x_i,f_1(x_1),...,\allowmathbreak f_i(x_1,...,x_i))=0$ as elements of $\Cal O(X_1\times\cdots
\times X_i)$. 
Then there exist $\tilde f_i(x_1,...,x_i)\in \Cal N(X_1\times\cdots\times X_i)^{l_i}$ close to $f_i(x_1,...,x_i)$ 
for $i=1,...,m$ such that $F_i(x_1,...,x_i,\tilde f_1(x_1),...,\tilde f_i(x_1,...,x_i))\allowmathbreak=0$ in 
$\Cal N(X_1\times\cdots\times X_i)$. 
\endproclaim

Let $M\supset X$ be a Nash manifold and a closed semialgebraic subset. 
Let $G(X)$ denote the germs on $X$ of semialgebraic subsets of $M$. 
Note that for $Z\in G(X)$, the topological closure $\overline Z$ is well-defined as an element of $G(X)$ since 
the topological closure of a semialgebraic set is semialgebraic. 
A {\it Nash set germ} on $X$ or in $G(X)$ is the germ on $X$ of a Nash subset of open semialgebraic 
neighborhood of $X$ in $M$. 
For an element $(\cdots)$ of $G(X)$, let $\overline{(\cdots)}^X$ or $(\cdots\overline)^X$ denote the smallest 
Nash set germ in $G(X)$ containing $(\cdots)$. 
We define by induction a sequence of elements $Z_1,Z_2,...$ of $G(X)$. 
Let $Z_1=\overline X^X$, and assume that $Z_1,...,Z_i$ are defined for some $i\,(>0)\in\N$. 
Then set 
$$
Z_{i+1}=[\overline{(Z_i-X)}\cap\overline{(Z_i\cap X)}\,\overline]^X. 
$$
We call $\{Z_i\}$ the {\it canonical Nash germ decomposition} of $X$. 
We see easily that $\{Z_i\}$ is a decreasing and finite sequence of Nash set germs in $G(X)$, and for each 
$i$ $X\cap Z_i-Z_{i+1}$ is the union of some connected components of $Z_i-Z_{i+1}$. 
Moreover, the canonical Nash germ decomposition is analytically invariant in the following sense. 
\remark{Remark 3 {\rm (Remark 3.3 in [F-S$_2$])}}
Let $M\supset X$ and $L\supset Y$ be Nash manifolds and closed semialgebraic subsets and $\phi\in\Cal O(
X,Y)$ an analytic diffeomorphism germ. 
Then $\phi$ carries the canonical global Nash germ decomposition of $X$ to the one of $Y$.
\endremark
\demo{Proof of theorem 1}
Proof contains two ideas. 
First we reduce the problem to the case where theorem 2 is applicable. 
Here we use essentially the fact that the image of a semialgebraic set under a semialgebraic map is semialgebraic. 
By this idea only theorem 4 is proved. 
Hence I recommend to read the proof of theorem 4 before the present proof. 
Secondly we use the sheaf theory and gather globally the local data by Cartan Theorem on Stein manifolds and the corresponding Nash case ([C-R-S] and [C-S]). \par
  Let $\tau_1\in\Cal N(X_1,Y_1)$ and $\tau_2\in\Cal N(X_2,Y_2)$ be Nash diffeomorphism germs such that $g\circ
\tau_1=\tau_2\circ f$. 
$$
\CD
G(X_2)@>\tau_2>>G(Y_2)@.\qquad\qquad M_2@>\hat\tau_2>>L_2\\
@A f AA@AA g A \qquad\qquad@A\hat f AA@AA\hat g A\\
G(X_1)@>\tau_1>>G(Y_1)@.\qquad\qquad M_1@>\hat\tau_1>>L_1
\endCD
$$
Let $\{Z_{X_1,i}\},\,\{Z_{X_2,i}\},\,\{Z_{Y_1,i}\}$and $\{Z_{Y_2,i}\}$ be the canonical Nash germ decompositions 
of $X_1,\,X_2,\,Y_1$ and $Y_2$, respectively. 
Then by remark 3
$$
\tau_j(Z_{X_j,i})=Z_{Y_j,i}\quad\text{for any}\ i\ \text{and}\ j. \tag 1
$$
Shrink $M_1,\,M_2,\,L_1$ and $L_2$ if necessary. 
Then we can assume that $Z_{X_j,i}$ and $Z_{Y_j,i}$ are the germs on $X_j$ and $Y_j$ of some Nash sets $\hat Z_{
X_j,i}$ in $M_j$ and $\hat Z_{Y_j,i}$ in $L_j$, respectively, $f$ and $g$ are the germs on $X_1$ and $Y_1$ of 
some Nash maps $\hat f:M_1\to M_2$ and $\hat g:L_1\to L_2$, respectively, and $\tau_j$ are the germs on $X_j$ of 
some analytic imbeddings $\hat\tau_j:M_j\to L_j,\,j=1,2$, with $\hat g\circ\hat\tau_1=\hat\tau_2\circ\hat f$ and 
$$
\hat\tau_j(\hat Z_{X_j,i})=\hat Z_{Y_j,i}\cap\Im\hat\tau_j. \tag$\hat1$
$$
Moreover, let $M_j$ and $L_j$ be contained and closed in $\R^n$ as Nash submanifolds, and $h_j\in N(R^n)$ 
with zero set $L_j$ (see [S$_1$]). 
Set $\hat F=\graph\hat f$ and $\hat G=\graph\hat g$, and let $F$ and $G$ denote their respective germs on $X_1
\times X_2$ and $\R^n\times\R^n$. 
(Hence $\hat G=G$.) 
Then for any analytic imbedding $\hat\tau'_j:M_j\to L_j,\,j=1,2$, $\hat g\circ\hat\tau'_1=\hat\tau'_2\circ\hat f$ 
if and only if $(\hat2)$ $\hat\tau'_1\times\hat\tau'_2(\hat F)\subset\hat G$. \par 
¡¡Let $\ _x$ denote the germ of a function or the stalk of a sheaf at a point $x$. 
For a Nash set $Z$ in $\R^n$, let $Z^\C$ denote its complexification---a complex analytic set in an open 
neighborhood of $\R^n$ in $\C^n$ containing $Z$ whose germ on $\R^n$ is the smallest. 
Let $\Cal I_{X_j,i},\,\Cal I_{Y_j,i},\,\Cal J_F$ and $\Cal J_G$ be the sheaves of $\Cal N_{M_j}$-,\,$\Cal N_{\R^n}$-,
\,$\Cal N_{M_1\times M_2}$- and $\Cal N_{\R^n\times\R^n}$-ideals of germs whose complexifications are vanishing on 
$\hat Z_{X_j,i}^\C,\,\hat Z_{Y_j,i}^\C,\,\hat F^\C$ and $\hat G^\C$, respectively. 
Then there exist global generators $\{\hat\alpha_{X_j,i,k}:k=1,...,m\}$ of $\Cal I_{X_j,i}$, $\{\hat\alpha_{Y_j,i,k}
:k=1,...,m\}$ of $\Cal I_{Y_j,i}$, $\{\hat\beta_{F,k}:k=1,...,m\}$ of $\Cal J_F$ and $\{\hat\beta_{G,k}:k=1,...,m
\}$ of $\Cal J_G$ by [C-R-S] and [C-S]. 
Remember the fact, which follows from AA Theorem, Nullstellensatz and the faithfully flatness of the completion of 
the local ring $\Cal O_{M_1x_1}$ over $\Cal O_{M_1x_1}$ for $x_1\in M_1$, that the ideal of $\Cal N_{M_1x_1}$ of 
germs whose complexifications are vanishing on the complexification of a Nash set germ at $x_1$ generates the ideal 
of $\Cal O_{M_1x_1}$ of germs whose complexifications are vanishing on the complexification of the Nash set germ. 
By the fact, $(\hat1)$ and $(\hat2)$ for $\hat\tau_1$ and $\hat\tau_2$, $\hat\alpha_{Y_j,i,k}\circ\hat\tau_j$ and 
$\hat\beta_{G,k}\circ\hat\tau_1\times\hat\tau_2$ are global cross-sections of $\Cal I_{X_j,i}\Cal O_{M_j}$ and 
$\Cal J_F\Cal O_{M_1\times M_2}$, respectively. 
Apply Cartan Theorem B on Stein manifolds to the homomorphisms 
$$
\gather
\Cal O_{M_j}^m\supset\Cal O_{M_jx_j}^m\ni(\xi_1,...,\xi_m)\to\sum_{k=1}^m\xi_k\hat\alpha_{X_j,i,kx_j}\in\Cal I_{X_j,
ix_j}\Cal O_{M_jx_j}\subset\Cal I_{X_j,i}\Cal O_{M_j},\\
\Cal O_{M_1\times M_2}^m\supset\Cal O_{M_1\times M_2(x_1,x_2)}^m\ni(\xi_1,...,\xi_m)\to\qquad\qquad\qquad\qquad
\qquad\qquad\qquad\\
\qquad\qquad\qquad\qquad\sum_{k=1}^m\xi_k\hat\beta_{F,k(x_1,x_2)}\in\Cal J_{F(x_1,x_2)}\Cal O_{M_1\times M_2(x_1,
x_2)}\subset\Cal J_F\Cal O_{M_1\times M_2}\\
\qquad\qquad\qquad\qquad\qquad\qquad\qquad\qquad\text{for}\ j=1,2\ \text{and}\ (x_1,x_2)\in M_1\times M_2. 
\endgather
$$ 
Then the induced homomorphisms $H^0(M_j,\Cal O^m_{M_j})\to H^0(M_j,\Cal I_{X_j,i}\Cal O_{M_j})$ and $H^0(M_1
\allowmathbreak\times M_2,\Cal O^m_{M_1\times M_2})\to H^0(M_1\times M_2,\Cal J_F\Cal O_{M_1\times M_2})$ are 
surjective. 
Hence there exist analytic functions $\hat\xi_{j,i,k,k'}$ on $M_j$ and $\hat\eta_{k,k'}$ on $M_1\times M_2$ such that 
$$
\gather 
\hat\alpha_{Y_j,i,k}\circ\hat\tau_j=\sum_{k'=1}^m\hat\xi_{j,i,k,k'}\hat\alpha_{X_j,i,k'}\quad\text{on}\ M_j,\\
\hat\beta_{G,k}\circ\hat\tau_1\times\hat\tau_2=\sum_{k'=1}^m\hat\eta_{k,k'}\hat\beta_{F,k'}\quad\text{on}\ M_1\times 
M_2. 
\endgather
$$\par
  Let $\alpha_{Y_j,i,k},\,\xi_{j,i,k,k'},\,\alpha_{X_j,i,k'},\,\beta_{G,k},\,\eta_{k,k'}$ and $\beta_{F,k}$ be the 
germs of $\hat\alpha_{Y_j,i,k}$ on $\R^n$, of $\hat\xi_{j,i,k,k'}$ on $X_j$, of $\hat\alpha_{X_j,i,k'}$ on $X_j$, 
of $\hat\beta_{G,k}$ on $\R^n\times\R^n$, of $\hat\eta_{k,k'}$ on $X_1\times X_2$ and of $\hat\beta_{F,k}$ on 
$X_1\times X_2$, respectively. 
Let $(z_1,z_2)\in\R^n\times\R^n$, $u=(u_{j,i,k,k'})\in\R^{2N}$ and $v=(v_{k,k'})\in\R^{N'}$, where $N=m^2\#\{Z_{X_1,
i}\}$ and $N'=m^2$. 
Consider the germ on $X_1\times X_2\times\R^n\times\R^n\times\R^{2N}\times\R^{N'}$ of the following Nash function on $M_1\times M_2\times\R^n\times\R^n\times\R^{2N}\times\R^{N'}$:
$$
\gather
K(x_1,x_2,z_1,z_2,u,v)=\sum_{j,i,k}[\alpha_{Y_j,i,k}(z_j)-\sum_{k'=1}^mu_{j,i,k,k'}\alpha_{X_j,i,k'}(x_j)]^2+\qquad
\qquad\qquad\\
\sum_k[\beta_{G,k}(z_1,z_2)-\sum_{k'=1}^mv_{k,k'}\beta_{F,k'}(x_1,x_2)]^2+h_1^2(z_1)+h^2_2
(z_2)\\
\qquad\qquad\qquad\qquad\text{for}\ (x_1,x_2,z_1,z_2,u,v)\in M_1\times M_2\times\R^n\times\R^n\times\R^{2N}\times\R^{N'}.
\endgather
$$
Then $(z_1,z_2,u_{1,i,k,k'},u_{2,i,k,k'},v_{k,k'})\!=\!(\tau_1(x_1),\tau_2(x_2),\xi_{1,i,k,k'}(x_1),\xi_{2,i,k,k'}(x_2),
\eta_{k,k'}(x_1,x_2))\allowmathbreak\in\Cal O(X_1)^n\times\Cal O(X_2)^n\times\Cal O(X_1)^N\times\Cal O(X_2)^N\times
\Cal O(X_1\times X_2)^{N'}$ is an analytic solution of the equation $K=0$. 
Hence by NGAA Theorem there exists a Nash solution $(\tilde\tau_1(x_1,x_2),\tilde\tau_2(x_2),\allowmathbreak\tilde\xi_{1,i,k,k'}
(x_1,x_2),\tilde\xi_{2,i,k,k'}(x_2),\tilde\eta_{k,k'}(x_1,x_2))\in\Cal N(X_1\times X_2)^n\times\Cal N(X_2)^n\times
\Cal N(X_1\times X_2)^N\times\Cal N(X_2)^N\times\Cal N(X_1\times X_2)^{N'}$ close to the analytic solution. 
Then 
$$
\gather
\alpha_{Y_1,i,k}(\tilde\tau_1(x_1,x_2))=\sum_{k'=1}^m\tilde\xi_{1,i,k,k'}(x_1,x_2)\alpha_{X_1,i,k'}(x_1),\tag 3\\
\alpha_{Y_2,i,k}(\tilde\tau_2(x_2))=\sum_{k'=1}^m\tilde\xi_{2,i,k,k'}(x_2)\alpha_{X_2,i,k'}(x_2),\tag 4\\
\beta_{G,k}(\tilde\tau_1(x_1,x_2),\tilde\tau_2(x_2))=\sum_{k'=1}^m\tilde\eta_{k,k'}(x_1,x_2)\beta_{F,k}(x_1,x_2),
\tag 5\\
(6)\,\qquad\qquad\quad h_1(\tilde\tau_1(x_1,x_2))=0,\qquad\quad(7)\qquad\qquad h_2(\tilde\tau_2(x_2))=0.\ \qquad
\qquad\qquad
\endgather
$$\par
  Set $\Tilde{\Tilde\tau}_1(x_1)=\tilde\tau_1(x_1,f(x_1))$ and $\Tilde{\Tilde\xi}_{1,i,k,k'}(x_1)=\tilde\xi_{1,i,k,
k'}(x_1,f(x_1))$. 
Then $\Tilde{\Tilde\tau}_1$ is close to $\tau_1$, and (3) and (6) imply, respectively, 
$$
\gather
\alpha_{Y_1,i,k}(\Tilde{\Tilde\tau}_1(x_1))=\sum_{k'=1}^m\Tilde{\Tilde\xi}_{1,i,k,k'}(x_1)\alpha_{X_1,i,k'}(x_1),
\tag $3'$\\
h_1(\Tilde{\Tilde\tau}_1(x_1))=0.
\tag$6'$
\endgather
$$
Hence $\Im\Tilde{\Tilde\tau}_1$ is contained in $L_1$ by $(6')$, each $\Tilde{\Tilde\tau}_1(Z_{X_1,i})$ is 
contained in $Z_{Y_1,i}$ by $(3')$, and $\Tilde{\Tilde\tau}_1$ is an imbedding since $\Tilde{\Tilde\tau}_1$ is close 
to the diffeomorphism germ $\tau_1$. 
Therefore, $\Tilde{\Tilde\tau}_1(Z_{X_1,i})=Z_{Y_1,i}$ and by the above-mentioned property of the canonical Nash 
germ decomposition, $\Tilde{\Tilde\tau}_1(X_1)=Y_1$. 
Hence $\Tilde{\Tilde\tau}_1$ is a Nash diffeomorphism germ in $\Cal N(X_1,Y_1)$. 
It follows also from (4) and (7) that $\tilde\tau_2$ is a Nash diffeomorphism germ in $\Cal N(X_2,Y_2)$. \par
  It remains to see $g\circ\Tilde{\Tilde\tau}_1=\tilde\tau_2\circ f$, i.e., $\Tilde{\Tilde\tau}_1\times\tilde\tau_2(F
)\subset G$. 
Shrink, moreover, $M_1$ and $M_2$ so that $(\tilde\tau_1,\tilde\tau_2)$ and $\tilde\eta_{k,k'}$ are the germs on $X_1
\times X_2$ of some Nash imbedding $(\check\tau_1,\check\tau_2):M_1\times M_2\to L_1\times L_2$ and of some Nash 
functions $\check\eta_{k,k'}$ on $M_1\times M_2$, respectively, with $\check\tau_2$ in only the variable $x_2\in M_2$ 
and 
$$
\hat\beta_{G,k}(\check\tau_1(x_1,x_2),\check\tau_2(x_2))=\sum_{k'=1}^m\check\eta_{k,k'}(x_1,x_2)\hat\beta_{F,k}(x_1,
x_2)\quad\text{on}\ M_1\times M_2, \tag$\check5$
$$
which follows from (5). 
Then $(\check5)$ implies that $(\check\tau_1,\check\tau_2)$ carries $\hat F$ into $\hat G$. 
Hence, if $(x_1,x_2)\in\hat F$ then $x_2=\hat f(x_1)$, $(\check\tau_1(x_1,x_2),\check\tau_2(x_2))\in\hat G$ and 
$(\check\tau_1(x_1,\hat f(x_1)),\check\tau_2(x_2))\in\hat G$. 
Therefore, $\Tilde{\Tilde\tau}_1\times\tilde\tau_2(F)\subset G$. 
\qed
\enddemo

\head \S 3. Equivalence of map germs at a point
\endhead
We naturally define {\it Nash, analytic} and $C^\infty$ {\it R-L} and {\it R equivalence} of two map germs 
$:(\R^n,0)\to(\R^m,0)$ by diffeomorphism germs $:(\R^n,0)\to(\R^n,0)$ and $:(\R^m,0)\to(\R^m,0)$. 
We say that $C^\infty$ map germs $f,g:(\R^n,0)\to(\R^m,0)$ are {\it formally R-L (R) equivalent} if there exist 
$C^\infty$ diffeomorphism germs $\tau_1:(\R^n,0)\to(\R^n,0)$ and $\tau_2:(\R^m,0)\to(\R^m,0)$ such that 
$T(g\circ\tau_1)=T(\tau_2\circ f)$ (and $\tau_2=\id$, respectively), where $T$ denotes the Taylor expansion at 0. 
Clearly formal R-L and R equivalence is weaker than respective $C^\infty$ one. 
As noted, $C^\infty$ R-L equivalence of two global Nash maps does not imply Nash R-L equivalence. 
In this section we consider whether $C^\infty$ R-L equivalence of two Nash or analytic map germs implies Nash or 
analytic R-L equivalence, respectively. \par
  The case of R equivalence is easy to see. 
2 $C^\infty$ R equivalent analytic map germs at 0 are formally R equivalent and hence analytically R equivalent 
by another AA Theorem [A$_1$], 
which says that a formal solution of a local analytic equation is approximated by an analytic solution. 
By this fact and AA Theorem in [A$_2$], two $C^\infty$ R equivalent Nash map germs at 0 are Nash R equivalent. \par
  The answer to first problem on Nash map germs is positive. 
\proclaim{Theorem 4} 
2 formally R-L equivalent Nash map germs at 0 are Nash R-L equivalent. 
\endproclaim

To prove theorem 4 we use the next NAA Theorem by Teissier. 
Let $x_1\in\R^{n_1},...,x_m\in\R^{n_m},\,y_1\in\R^{l_1},...,y_m\in\R^{l_m}$. 
Let $\R[[\cdots]]$ and $\R_{\text{alg}}[[\cdots]]$ denote the rings of formal power series and Nash function germs 
at 0, respectively, with the $\frak p$-topology, where $\frak p$ denotes the maximal ideals. 
An element $\tau$ of $\R[[x_1]]^{n_1}$ is called {\it invertible} if $\tau(0)=0$ and if there exists $\pi\in\R[[x_1]]^{n_1}$ such that $\pi(0)=0$ and $\tau\circ\pi(x_1)=\pi\circ\tau(x_1)=x_1$, i.e., $\tau$ is the Taylor expansion at 0 of a $C^\infty$ diffeomorphism germ $:(\R^{n_1},0)\to(\R^{n_1},0)$. 

\proclaim{Theorem 5} {\rm(See [Sp].)} 
Let $F_i\in\R_{\text{alg}}[[x_1,...,x_i,y_1,...,y_i]]$ for $i=1,...,m$ and $f_i\in\R[[x_1,...,x_i]]^{l_i}$ such that $F_i(x_1,...,x_i,f_1(x_1),...,f_i(x_1,...,x_i))=0$. 
Then $f_i$ are approximated by $\tilde f_i\in\R_{\text{alg}}[[x_1,...,x_i]]^{l_i}$ so that $F_i(x_1,...,x_i,\tilde f
_1(x_1),...,\tilde f_i(x_1,...,x_i))=0$. 
\endproclaim
\demo{Proof of theorem 4}
Let $x=(x_1,...,x_n),\,u=(u_1,...,u_n)\in\R^n$ and $y=(y_1,...,y_m),\,\allowmathbreak v=(v_1,...,v_m)\in\R^m$. 
Let $y=f(x)=(f_1(x),...,f_m(x)),\,v=g(u)=(g_1(u),...,\allowmathbreak g_m(u)):(\R^n,0)\to(\R^m,0)$ be Nash map germs 
and $u=\tau(x)=(\tau_1(x),...,\tau_n(x))\in\R[[x_1,...,x_n]]^n$ and $v=\pi(y)=(\pi_1(y),...,\pi_m(y))\in\R[[y_1,...,
y_m]]^m$ invertible elements such that $\pi\circ f=g\circ\tau$. 
$$
\CD
y@>\pi>>v\\
@A f AA@AA g A\\
x@>\tau>>u
\endCD
$$
Then $\pi_i\circ f=g_i\circ\tau$ for each $i=1,...,m$. 
Hence $\pi_i(y)-g_i\circ\tau(x)=0$ if $y_j=f_j(x),\,j=1,...,m$. 
Regard $\pi_i(y)-g_i\circ\tau(x)$ and $y_j-f_j(x)$ as elements of $\R[[x,y]]$. 
Then $\{x_1,...,x_n,y_1-f_1(x),...,y_m-f_m(x)\}$ is a basis of $\R[[x,y]]$ and each $\pi_i(y)-g_i\circ\tau(x)$ 
is described as $\sum_{\alpha'\in\N^n,\alpha''\in\N^m}a_{\alpha',\alpha''}x^{\alpha'}(y-f(x))^{\alpha''},\ a_{
\alpha',\alpha''}\in\R$, where $\alpha'=(\alpha'_1,...,\alpha'_n)\in\N^n,\,\alpha''=(\alpha''_1,...,\alpha''_m)\in\N^m,\,x^
{\alpha'}=\prod_{i=1}^nx^{\alpha'_i}$ and $(y-f(x))^{\alpha''}=\prod_{j=1}^m(y_j-f_j(x))^{\alpha''_j}$. 
Such a description is unique. 
Hence $a_{\alpha',\alpha''}=0$ for $(\alpha',\alpha'')$ with $\alpha''=0$. 
Therefore, we have $\beta_{i,j}\in\R[[x,y]],\,i,j=1,...,m$, such that 
$$
\pi_i(y)-g_i\circ\tau(x)=\sum_{j=1}^m\beta_{i,j}(x,y)(y_j-f_j(x)),\ i=1,...,m.
$$\par
  Let $w=(w_{i,j})\in\R^{m^2}$ and consider the following power series of class Nash 
$$
F(x,y,u,v,w)=\sum_{i=1}^m[v_i-g_i(u)-\sum_{j=1}^mw_{i,j}(y_j-f_j(x))]^2.
$$
Then $\{u=\tau(x),\,v=\pi(y),\,w=(\beta_{i,j}(x,y))\}$ is a solution of formal power series of the equation 
$F(x,y,u,v,w)=0$. 
Hence by NAA Theorem we have a solution of power series of class Nash $\{u=\tilde\tau(x,y)=(\tilde\tau_1(x,y),...,
\tilde\tau_n(x,y)),\,v=\tilde\pi(y)=(\tilde\pi_1(y),...,\tilde\pi_m(y)),\,w=(\tilde\beta_{i,j}(x,y))\}$ close 
to $\{\tau(x),\,\pi(y),\,(\beta_{i,j}(x,y))\}$. 
Then 
$$
\gather
\tilde\pi_i(y)-g_i\circ\tilde\tau(x,y)=\sum_{j=1}^m\tilde\beta_{i,j}(x,y)(y_j-f_j(x)),\ i=1,...,m,\\
\tilde\pi_i\circ f(x)-g_i\circ\tilde\tau(x,f(x))=0.\tag"hence"
\endgather
$$
Set $\Tilde{\Tilde\tau}(x)\tilde\tau(x,f(x))$. 
Then $\tilde\pi\circ f(x)-g\circ\Tilde{\Tilde\tau}(x)=0$ and $\Tilde{\Tilde\tau}$ is close to $\tau$ and hence 
invertible. 
Thus $f$ and $g$ are Nash R-L equivalent and theorem 4 is proved. 
\qed
\enddemo

The answer to the second problem on analytic map germs is negative. 

\example{Example}
Let $f,g:(\R^2,0)\to(\R^4,0)$ be analytic map germs defined by 
$$
\gather
f(x_1,x_2)=(x_1,x_1x_2,x_1x_2e^{x_2},0),\\
g(x_1,x_2)=(x_1,x_1x_2,x_1x_2e^{x_2},\sum_{k=1}^\infty\sum_{i=0}^\infty\frac{k!}{(k+i)!}x_1^kx_2^{k+i+1}).
\endgather
$$
Then $f$ and $g$ are $C^\infty$ L equivalent but not $C^\omega$ R-L equivalent, where $C^\infty$ 
{\it L equivalence} is defined by only a diffeomorphism germ of the target space. 
\endexample
This example comes from a counter-example to analytic NAA Theorem (theorem 5 with convergent power series $F_i$) by [G]. 
The author already claimed in [S$_3$] that $f$ and $g$ are $C^\infty$ L equivalent and not $C^\omega$ R-L 
equivalent, but there is a gap in the proof of $C^\infty$ L equivalence. 
We correct here it. 

\demo{Proof of $C^\infty$ L equivalence}
Let $y=(y_1,..,y_4)\in\R^4$, and define $\pi\in\R[[y]]^4$ by
$$
\pi(y)=(y_1,y_2,y_3,\pi_4(y))\quad\text{and}\ \ \pi_4(y)=y_4-\sum_{k=1}^\infty(k!y_1^{k-1}y_3-\sum_{i=1}^k
\frac{k!}{(i-1)!}y_1^{k-i}y_2^i). 
$$
Then $\pi$ is invertible and 
$$
\gather
\pi\circ g(x)=(x_1,x_1x_2,x_1x_2e^{x_2},\qquad\qquad\qquad\qquad\qquad\qquad\qquad\qquad\qquad\qquad\\
\sum_{k=1}^\infty\sum_{i=0}^\infty\frac{k!}{(k+i)!}x_1^kx_2^{k+i+1}-\sum_{k=1}^\infty[k!x_1^{k-1}(x_1x_2
\sum_{i=0}^\infty\frac{x_2^i}{i!})-\sum_{i=1}^k\frac{k!}{(i-1)!}x_1^{k-i}(x_1x_2)^i])=\\
\qquad\qquad\qquad\qquad\qquad\qquad\qquad\qquad\qquad\qquad(x_1,x_1x_2,x_1x_2e^{x_2},0)=f(x). 
\endgather
$$
Hence $f$ and $g$ are formally L equivalent. \par
  Let $\tilde\pi,\tau:(\R^4,0)\to(\R^4,0)$ be $C^\infty$ diffeomorphism germs of the form $\tilde\pi(y)=(y_1,
  y_2,y_3,\allowmathbreak\tilde\pi_4(y))$ and $\tau(y)=(y_1,y_2,y_3,y_4+\tau_4(y_1,y_2))$ and such that $T\tilde\pi_4=\pi_4$. 
Then 
$$
\gather
\tilde\pi\circ g(x)=(x_1,x_1x_2,x_1x_2e^{x_2},\tilde\pi_4\circ g(x))\quad\text{and}\\
\tau\circ f(x)=(x_1,x_1x_2,x_1x_2e^{x_2},\tau_4(x_1,x_1x_2)). 
\endgather
$$
Hence it suffices to find $C^\infty$ function germs $\tilde\pi_4(y)$ and $\tau_4(y_1,y_2)$ such that $T\tilde
\pi_4=\pi_4$ and $\tau_4(x_1,x_1x_2)=\tilde\pi_4\circ g(x)$. \par
  First we define $\tilde\pi_4$. 
Let $\phi$ be a $C^\infty$ function on $\R$ such that $\phi=0$ outside of $[-1,\,1]$ and $\phi=1$ on $[-1/2,\,1/2]$. 
Set 
$$
\gather\tilde\pi_4(y)=y_4-\sum_{k=1}^\infty\phi(ky_1)(k!y_1^{k-1}y_3-\sum_{i=1}^k\frac{k!}{(i-1)!}y_1^{k-i}y_2^i)
\qquad\qquad\qquad\qquad\\
\qquad\qquad\qquad\qquad\text{for}\ y=(y_1,..,y_4)\in\R\times(-1/2,\,1/2)\times\R^2,
\endgather
$$
which is a well-defined $C^\infty$ function for the following reason. 
Let $\alpha=(\alpha_1,\alpha_2,\alpha_3)\in\N^3$, $0<l_1<l_2\in\N$, and $K$ a compact subset of $\R\times
(-1/2,\,1/2)\times\R$. 
Then we need to see that the restrictions to $K$ of the following functions uniformly converge to 0 as $l_1,l_2
\to\infty$\,: 
$$
\gather
\xi_{1,l_1,l_2}(y_1,y_3)=\sum_{k=l_1}^{l_2}\phi^{(\alpha_1)}(ky_1)k^{\alpha_1}k!(k-1)(k-2)\cdots(k-\alpha_2)y_1
^{k-\alpha_2-1}y_3,\\
\xi_{2,l_1,l_2}(y_1)=\sum_{k=l_1}^{l_2}\phi^{(\alpha_1)}(ky_1)k^{\alpha_1}k!(k-1)(k-2)\cdots(k-\alpha_2)y_1^{k-
\alpha_2-1},\\
\xi_{3,l_1,l_2}(y_1,y_2)=\sum_{k=l_1}^{l_2}\sum_{i=\alpha'_3}^{k-\alpha_2}\phi^{(\alpha_1)}(ky_1)\times\qquad
\qquad\qquad\qquad\qquad\qquad\qquad\qquad\\\frac{k^{\alpha_1}k!(k-i)(k-i-1)\cdots(k-i-\alpha_2+1)i(i-1)\cdots
(i-\alpha_3+1)}{(i-1)!}y_1^{k-i-\alpha_2}y_2^{i-\alpha_3}, 
\endgather
$$
where $\alpha'_3=\max\{1,\alpha_3\}$. 
If $\phi^{(\alpha_1)}(ky_1)\not=0$ then $|y_1|\le k^{-1}$. 
Hence 
$$
|\xi_{1,l_1,l_2}|_K|\le\sum_{k=l_1}^{l_2}ck!k^{-k+c'}\quad\text{for some constants}\ c,\,c'\in\N. 
$$
Choose $l_1$ so large that $(c'+3)!\le l_1$. 
Then $(c'+3)!k^{-k+c'}\le k^{-k+c'+1}$ for $k\ge l_1$. 
Hence 
$$
\gather
\sum_{k=l_1}^{l_2}ck!k^{-k+c'}\le\sum_{k=l_1}^{l_2}c\oversetbrace k-c'-3\to{k(k-1)\cdots(c'+4)}k^{-k+c'+3}k^{-2}
\le\sum_{k=l_1}^{l_2}ck^{-2}\to0\\
\qquad\qquad\qquad\qquad\qquad\qquad\qquad\qquad\text{as}\ l_1,\,l_2\to\infty. 
\endgather
$$
Theses arguments show also that $|\xi_{2,l_1,l_2}(y_1)|\to0$ as $l_1,\,l_2\to\infty$. \par
  Consider $\xi_{3,l_1,l_2}$. 
Since $|y_2|<2^{-1}$, we have 
$$
\gather
|\xi_{3,l_1,l_2}(y_1,y_2)|\le\sum_{k=l_1}^{l_2}\sum_{i=\alpha'_3}^{k-\alpha_2}ck(k-1)\cdots ik^{-k+i+c'}2^{-i}\\
\le\sum_{k=l_1}^{l_2}\max_{\alpha'_3\le i\le k-\alpha_2}(k-\alpha_2-
\alpha'_3+1)ck(k-1)\cdots ik^{-k+i+c'}2^{-i}\\
\le\sum_{k=l_1}^{l_2}\max_{\alpha'_3\le i\le k-\alpha_2}\frac{ck(k-1)\cdots ik^{-k+i+c'+3}2^{-i}}{k^2}
\endgather
$$
for $l_1\ge\alpha_2+\alpha'_3$ and for some constants $c,c'\in\N$. 
Hence it suffices to see $k(k-1)\cdots ik^{-k+i+c'+3}2^{-i}$ is bounded for all $i$ and $k$ with $0<i\le k$. 
For each $i>0\in\N$, let $a_i$ be the largest natural number such that $(i+a_i)(i+a_i-1)\cdots i\le2^i$. 
Then $a_i\to\infty$ as $i\to\infty$. 
Hence there are only a finite number of $i$'s with $a_i<c'+3$. 
For such $i$, the boundedness of $k(k-1)\cdots ik^{-k+i+c'+3}2^{-i}$ follows from the above arguments on 
$\xi_{1,l_1,l_2}$. 
Assume that $a_i\ge c'+3$. 
There are two possible cases to consider\,$:k\ge i+a_i$ or $k<i+a_i$. 
If $k\ge i+a_i$, then 
$$
k(k-1)\cdots ik^{-k+i+c'+3}2^{-i}\le\oversetbrace k-i-a_i\to{k(k-1)\cdots(i+a_i+1)}k^{-k+i+a_i}\le1. 
$$
If $k<i+a_i$, then 
$$
\gather
k(k-1)\cdots ik^{-k+i+c'+3}2^{-i}=\frac{(i+a_i)(i+a_i-1)\cdots ik^{-k+i+c'+3}2^{-i}}{(i+a_i)(i+a_i-1)
\cdots(k+1)}\qquad\\
\qquad\qquad\qquad\qquad\qquad\qquad\le\frac{k^{-k+i+a_i}}{\undersetbrace -k+i+a_i\to{(i+a_i)(i+a_i-1)
\cdots(k+1)}}\le1. 
\endgather
$$\par
  Thus $\tilde\pi_4$ is a $C^\infty$ function on $\R\times(-1/2,\,1/2)\times\R^2$. 
It is clear that $T\tilde\pi_4(y)=\pi_4(y)$. \par
  Next we find $\tau_4(y_1,y_2)$. 
By calculations we have 
$$
\multline
\tilde\pi_4\circ g(x)=\\
\sum_{k=1}^\infty\sum_{i=0}^\infty\frac{k!}{(k+i)!}x_1^kx_2^{k+i+1}-\sum_{k=1}^\infty\phi(kx_1)[k!x_1^{k-1}
(x_1x_2e^{x_2})-\sum_{i=1}^k\frac{k!}{(i-1)!}x_1^{k-i}(x_1x_2)^i]\\
=\sum_{k=1}^\infty(1-\phi(kx_1))\sum_{i=0}^\infty\frac{k!}{(k+i)!}x_1^kx_2^{k+i+1}. 
\endmultline
$$
Hence the required equation $\tau_4(x_1,x_1x_2)=\tilde\pi_4\circ g(x)$ becomes 
$$
\tau_4(x_1,x_2)=\sum_{k=1}^\infty(1-\phi(kx_1))\sum_{i=0}^\infty\frac{k!x_2^{k+i+1}}{(k+i)!x_1^{i+1}}\,. 
$$
Define $\tau_4(x_1,x_2)$ for $x_1\not=0$ by this equality, and set $\tau_4(x_1,x_2)=0$ for $x_1=0$. 
Then $\tau_4$ is a $C^\infty$ function on $\R\times(-1/2,\,1/2)$. 
For that we only need to show that $\sum_{i=0}^N\frac{k!(2k)^{i+1+c}}{(k+i)!2^{k+i+1-c'}}$ converges as 
$N\to\infty$ and $\sum_{k=l_1}^{l_2}\sum_{i=0}^\infty\frac{k!(2k)^{i+1+c}}{(k+i)!2^{k+i+1-c'}}$ converges to 
0 as $l_1,l_2\to\infty$ for any given $c,c'\in\N$ by the same reason as before because $1-\phi(kx_1)=0$ for 
$x_1$ in $[-(2k)^{-1},\,(2k)^{-1}]$. 
Set $b=\max_kk^{5+c}2^{-k+c+c'}$. 
Let $k>1$. 
If $i>1$, then we have 
$$
\gather
\frac{k!(2k)^{i+1+c}}{(k+i)!2^{k+i+1-c'}}=\frac{k^{i-2}}{k^2(k+i)(k+i-1)\cdots(k+1)}\frac{k^{5+c}}{2^{k-c-c'}}
\qquad\qquad\qquad\qquad\\
\le\frac{bk^{i-2}}{k^2(k+i)(k+i-1)\cdots(k+1)}=\frac{b}{k^2(k+i)(k+i-1)}\frac{k^{i-2}}{\undersetbrace i-2\to{
(k+i-2)\cdots(k+1)}}\\
\qquad\qquad\qquad\qquad\qquad\qquad\qquad\qquad\le\frac{b}{k^2(k+i)(k+i-1)}<\frac{b}{k^2(i+1)^2}. 
\endgather
$$
If $i=0$ or $=1$, we see in the same way that 
$$
\frac{k!(2k)^{i+1+c}}{(k+i)!2^{k+i+1-c'}}<\frac{b}{k^4}. 
$$
Hence the two convergence properties follow. 
Thus $\tau_4$ is a $C^\infty$ function on $\R\times(-1/2,\,1/2)$. \linebreak
\qed
\enddemo
\remark{Problem}
If two analytic map germs $:(\R^n,0)\to(\R^m,0)$ are formally R-L equivalent then they are $C^\infty$ R-L equivalent. 
\endremark 
A partial answer is Fact 1.7 in [S$_3$].


\Refs
\widestnumber\key{C-R-S}
\ref 
\key A$_1$\by M.~Artin\paper On the solutions of analytic equations\jour Invent\. Math\.,\vol 5\page 277--291
\yr 1968\endref 
\ref 
\key A$_2$\bysame\paper Algebraic approximation of structures over complete local rings\jour Publ\. Math\. 
IHES,\vol 36\page 23--57\yr 1969\endref 
\ref
\key C-R-S\by M.~Coste, J.M.~Ruiz and M.~Shiota\paper Approximation in compact Nash manifolds
\jour Amer\. J\. Math\.,\vol 117\page 905--927\yr 1995\endref 
\ref
\key C-S\by M.~Coste and M.~Shiota\paper Nash functions on noncompact Nash manifolds\jour Ann\. 
Sci\. \'Ec\. Norm\. Sup\.,\vol 33\pages 139--149\yr 2000\endref
\ref
\key F-S$_1$\by G.~Fichou and M.~Shiota\paper Analytic equivalence of normal crossing functions on a real analytic manifold\jour preprint, arxiv:0811.4569
\vol \yr \pages \endref 
\ref
\key F-S$_2$\bysame \paper Nash approximation of an analytic desingularization\jour preprint, arxiv:0812.2335
\vol \yr \pages \endref 
\ref 
\key G\by A.~M.~Gabrielov\paper Formal relations between analytic functions\jour Functional Anal\. Appl\.,
\vol 5\yr 1971\pages 318-319\endref
\ref
\key S$_1$\by M.~Shiota\paper Nash manifolds\jour Lecture Notes in Math\.,\vol 1269\publ Springer\yr 1987
\endref
\ref
\key S$_2$\bysame\book Geometry of subanalytic and semialgebraic sets\publ Birkhauser\yr 1997
\endref
\ref
\key S$_3$\bysame\paper Relations between equivalence relations of maps and functions, Real analytic 
and algebraic singularities\jour Pitman Research Notes in Math.,\vol 381\yr 1998\pages 114--144\endref 
\ref
\key Sp\by M.~Spivakovsi\paper A new proof of D. Popescu's theorem on smoothing of ring homomorphisms
\jour J\. Amer\. Math\. Soc\.,\vol 12\yr 1999\pages 381--444\endref
\ref
\key T$_1$\by R.~Thom\paper La stabilit\'e topologique des applications polynomiales\jour Enseignement math\.,
\vol 8\yr 1962\pages 24--33\endref 
\ref
\key T$_2$\bysame \paper Ensembles et morphismes stratifi\'es\jour Bull\. AMS\.,
\vol 75\yr 1969\pages 240--284\endref 
\endRefs
\enddocument